\documentclass{article}
\usepackage{indentfirst,latexsym,bm,amsmath,amsthm,amsfonts}
\begin{document}
\title{Full graph and complete Lie algebra}

\author{BinYong HSIE
      \\{\small{\{ Department of Mathematics, PeKing University, BeiJing, 100871, P.R. China
      \}}}
      \\{\small{E-mail:xiebingyong1@sina.com.cn}}}

\date{2003.10}
\maketitle{}
\begin{abstract}
We introduce a concept: $d$-complete, and show that a Lie algebra
is $d$-complete if and only if its full graph is complete.
\end{abstract}


MR(2000) Subject Classification: 17B45, 17B55, 17B56

Let $\mathcal{G}$ be a Lie algebra. We will introduce two
concepts: $d$-derivation and $d$-center. Let $Der(\mathcal{G})$ be
the derivation algebra of $\mathcal{G}$.

A $d$-derivation of $\mathcal{G}$, means a linear operator $L$
from $Der(\mathcal{G})$ to $\mathcal{G}$ such that
$$L([D_{1},D_{2}])=D_{1}(L(D_{2}))-D_{2}(L(D_{1}))\;\;\forall D_{1},D_{2}\in Der(\mathcal{G}) \eqno(1).$$
We write $\mathcal{D}$ for the set of $d$-derivations.

The $d$-center of $\mathcal{G}$, means the set
$$C_{d}(\mathcal{G})\doteq \{g\in \mathcal{G}:D(g)=0\;\forall D\in Der(\mathcal{G})\}\eqno(2).$$

For any $g\in \mathcal{G}$, we can define a $d$-derivation $L_{g}$
of $\mathcal{G}$:
$$L_{g}(D)=-D(g)\;\;\forall D\in Der(\mathcal{G})\eqno(3).$$
Such a $d$-derivation is called inner $d$-derivation.

For any two $d$-derivations $L_{1}$,$L_{2}$, define
$[L_{1},L_{2}]$ by
$$[L_{1},L_{2}](D)=L_{1}(ad(L_{2}(D)))-L_{2}(ad(L_{1}(D)))\;\forall D\in Der(\mathcal{G})\eqno(4),$$
then $[L_{1},L_{2}]$ is a $d$-derivation. And it is easy to check
that ($\mathcal{D}$,[,]) is a Lie algebra. If $C_{d}(\mathcal{G})$
is trivial, $\mathcal{G}$ is a Lie subalgebra of $\mathcal{D}$.

If $C_{d}(\mathcal{G})$ is trivial and all $d$-derivations are
inner, $\mathcal{G}$ is called $d$-complete. In this case, there
is a isomorphic from $\mathcal{G}$ to $\mathcal{D}$: $g\mapsto
L_{g}$.

We have a homomorphic from $Der(\mathcal{G})$ to the derivation
algebra $Der(\mathcal{D})$ of $\mathcal{D}$ by setting:
$$D(L)\doteq D\circ L- L\circ ad(D)\;\forall L\in \mathcal{D},D\in
Der(\mathcal{G})\eqno(5).$$ We can now define semidirect product
$\mathcal{H}=Der(\mathcal{G})\times_{t}\mathcal{D}$, a Lie
algebra, by
$$[(D_{1},L_{1}),(D_{2},L_{2})]=([D_{1},D_{2}],[L_{1},L_{2}]+D_{1}(L_{2})-D_{2}(L_{1}))\eqno(6).$$

We can make $\mathcal{H}$ act on
$\mathrm{C}(\mathcal{G})=Der(\mathcal{G})\times_{t}\mathcal{D}$ as
derivations:
$$\mathrm{D}_{(D,L)}(D_{1},g)=([D,D_{1}],D(g)+L(ad(g))+L(D_{1}))\eqno(7).$$

In fact we have:

{\bf Theorem 1. $\mathcal{H}$ is the derivation algebra
$Der(\mathrm{C}(\mathcal{G}))$ of the full graph
$\mathrm{C}(\mathcal{G})$ of $\mathcal{G}$.}

\begin{proof}A direct calculation shows that

\textrm{1}.$\mathrm{D}_{(D,L)}$ is a derivation,

\textrm{2}.
$[\mathrm{D}_{(D_{1},L_{1})},\mathrm{D}_{(D_{2},L_{2})}]=\mathrm{D}_{[(D_{1},L_{1}),(D_{2},L_{2})]}$,

\textrm{3}. $\mathrm{D}_{(D,L)}=0$ if and only if $D=0$ and $L=0$.

So $\mathcal{H}$ is a Lie subalgebra of
$Der(\mathrm{C}(\mathcal{G}))$.

Now , we will show that every derivation $\mathrm{D}$ of
$Der(\mathrm{C}(\mathcal{G}))$ is in $\mathcal{H}$. We write
$Der(\mathcal{G})$ for $(Der(\mathcal{G}),0)$ and $\mathcal{G}$
for $(0,\mathcal{G})$. Write $\mathrm{D}=(T_{1},T_{2})$. First we
show that $T_{1}|\mathcal{G}=0$.

By
$$\mathrm{D}([(0,g_{1}),(D,g_{2})])=[\mathrm{D}(0,g_{1}),(D,g_{2})]+[(0,g_{1}),\mathrm{D}(D,g_{2})]\eqno(8),$$
we obtain
$$-T_{1}(D(g_{1}))+T_{1}([g_{1},g_{2}])=[T_{1}(g_{1}),D]\eqno(9).$$
Formula (9) shows that
$$T_{1}|[\mathcal{G},\mathcal{G}]=0\eqno(10)$$
and that
$$-T_{1}(D(g))=[T_{1}(g),D]\eqno(11).$$
Therefore
$$[T_{1}(g), D]=0\;\forall D\in ad(\mathcal{G}).$$
So $T_{1}(g)=0$ and $T_{1}|\mathcal{G}=0$.

We have
$$T_{2}([D_{1},D_{2}])=D_{1}(T_{2}(D_{2}))-D_{2}(T_{2}(D_{1})).$$
So $T_{2}|Der(\mathcal{G})=L\in \mathcal{D}$.

By setting $\widetilde{\mathrm{D}}=\mathrm{D}-\mathrm{D}_{(0,L)}$,
we may assume that $T_{2}|Der(\mathcal{G})=0$. We know that
$T_{1}|Der(\mathcal{G})$ is a derivation of $Der(\mathcal{G})$.

Since $T_{1}|\mathcal{G}=0$, $T_{2}|\mathcal{G}$ is a derivation
$D^{'}$ of $\mathcal{G}$. Now
$$\mathrm{D}[(D,0),(0,g)]=[(D,0),\mathrm{D}(0,g)]+[\mathrm{D}(D,0),(0,g)],$$
shows that
$$T_{2}(D(g))=D(T_{2}(g))+(T_{1}(D))(g),$$
that is
$$T_{1}(D)=[D^{'},D].$$
So $\mathrm{D}=\mathrm{D}_{(D^{'},0)}$. We have now showed that
$\mathrm{D}\in \mathcal{H}$.
\end{proof}

We have a easy lemma:

{\bf Lemma. The center of $\mathrm{C}(\mathcal{G})$ is
$(0,C_{d}(\mathcal{G}))$}

We have the following corollary of theorem 1 and the above lemma:

{\bf Theorem 2. A Lie algebra is $d$-complete if and only if its
full graph is complete.}

{\bf References:}

{\small

[1] Hsie B.Y., On complete Lie algebra, Science in China, in
appear.

[2] Meng D.J., Complete lie algebra(in Chinese), BeiJing: Press of
science,2001.
 }

{\large D}{\small EPARTMENT OF} {\large M}{\small ATHEMATICS,}
{\large P}{\small EKING} {\large U}{\small NIVERSITY,} {\large
B}{\small EIJING,} {\large P.R.  C}{\small HINA.}

E-mail address: xiebingyong1@sina.com.cn

\end{document}